\documentclass[a4paper,11pt]{amsart}
\usepackage{graphicx}
\usepackage{mathptmx}
\usepackage{amsmath}
\usepackage{amssymb}
\usepackage{enumitem}
\usepackage{xcolor}
\usepackage{xparse}
\NewDocumentCommand{\eulerian}{omm}
 {%
  \genfrac<>{0pt}{}{#2}{#3}%
  \IfValueT{#1}{_{\!#1}}%
 }

 \newmuskip\pFqmuskip

\newcommand*\pFq[6][8]{%
  \begingroup 
  \pFqmuskip=#1mu\relax
  \mathchardef\normalcomma=\mathcode`,
  \mathcode`\,=\string"8000
  \begingroup\lccode`\~=`\,
  \lowercase{\endgroup\let~}\pFqcomma
  {}_{#2}F_{#3}{\left(\genfrac..{0pt}{}{#4}{#5}\bigg|#6\right)}%
  \endgroup
}
\newcommand{\pFqcomma}{{\normalcomma}\mskip\pFqmuskip}

\newtheorem{theorem}{Theorem}
\newtheorem{lemma}[theorem]{Lemma}
\newtheorem{corollary}[theorem]{Corollary}

\begin{document}

\title[Some identities on degenerate $r$-Stirling numbers via boson operators]{Some identities on degenerate $r$-Stirling numbers via boson operators}

\author{Taekyun  Kim}
\address{Department of Mathematics, Kwangwoon University, Seoul 139-701, Republic of Korea}
\email{tkkim@kw.ac.kr}

\author{DAE SAN KIM}
\address{Department of Mathematics, Sogang University, Seoul 121-742, Republic of Korea}
\email{dskim@sogang.ac.kr}


\begin{abstract}
Broder introduced the $r$-Stirling numbers of the first kind and of the second kind which enumerate restricted permutations and respectively restricted partitions, the restriction being that the first $r$ elements must be in distinct cycles and respectively in distinct subsets. Kim-Kim-Lee-Park constructed the degenerate $r$-Stirling numbers of both kinds as degenerate versions of them. The aim of this paper is to derive some identities and recurrence relations for the degenerate $r$-Stirling numbers of the first kind and of the second kind via boson operators. In particular, we obtain the normal ordering of a degenerate integral power of the number operator multiplied by an integral power of the creation boson operator in terms of boson operators where the degenerate $r$-Stirling numbers of the second kind appear as the coefficients.
\end{abstract}

\maketitle

\section{Introduction}
The $r$-Stirling numbers of the first kind and of the second kind enumerate restricted permutations and respectively restricted partitions, the restriction being that the first $r$ elements must be in distinct cycles and respectively in distinct subsets (see [2]). \par 
Carlitz initiated a study of degenerate versions of some special numbers and polynomials (see [3]), where the degenerate Bernoulli and Euler numbers were investigated. It is remarkable that in recent years intensive studies have been done for degenerate version of quite a few special polynomials and numbers and have yielded many interesting results (see [7-15]). \par 
The normal ordering of an integral power of the number operator $a^{+}a$ in terms of boson operators $a$ and $a^{+}$ can be written in the form
\begin{equation*}
(a^{+}a)^{k}=\sum_{l=0}^{k}S_{2}(k,l)(a^{+})^{l}a^{l},
\end{equation*}
where $S_{2}(k,l)$ are the Stirling numbers of the second kind (see \eqref{2}).
The normal ordering of the degenerate $k$th power of the number operator $a^{+}a$, namely $(a^{+}a)_{k,\lambda}$, in terms of boson operators $a$ and $a^{+}$ can be written in the form
\begin{equation*}
(a^{+}a)_{k,\lambda}=\sum_{l=0}^{k}S_{2,\lambda}(k,l)(a^{+})^{l}a^{l},\label{1}
\end{equation*}
where $S_{2,\lambda}(k,l)$ are the degenerate Stirling numbers of the second kind (see \eqref{5}) and the generalized falling factorials $(x)_{n,\lambda}$ are given by \eqref{3}.\par
The aim of this paper is to derive some identities and recurrence relations for the degenerate $r$-Stirling numbers of the first kind and of the second kind via boson operators. In particular, we obtain the normal ordering of a degenerate integral power of the number operator multiplied by an integral power of the creation boson operator in trems of boson operators where the degenerate $r$-Stirling numbers of the second kind appear as the coefficients. In addition, we derive recurrence relations for the degenerate $r$-Stirling numbers of both kinds from certain normal orderings which generalize aforementioned ones.\par
The outline of this paper is as follows. In Section 1, we recall the Stirling numbers of both kinds, the unsigned Stirling numbers of the first kind, the degenerate Stirling numbers of both kinds, and the unsigned degenerate Stirling numbers of the first kind. We remind the reader of the normal ordering of an integral power of the number operator in terms of boson operators and its degenerate version, namely the normal ordering of a degenerate integral power of the number operator. We also recall the degenerate $r$-Stirling numbers of both kinds and the unsigned degenerate $r$-Stirling numbers of the first kind. Section 2 is the main result of this paper. Let $D=\frac{d}{dx}$. In Theorem 1, we show that the degenerate differential operator $(xD)_{m,\lambda}x^{r}$ applied to the formal power series $f(x)$ can be expressed as a linear combination of $x^{l+r}D^{l}f(x)$ with the coefficients given by the degenarate $r$-Stirling numbers of the second kind. Thereby we obtain the normal ordering of $(a^{+}a)_{n,\lambda}(a^{+})^{r}$ in terms of boson operators (see Corollary 2):
\begin{displaymath}
(a^{+}a)_{n,\lambda}(a^{+})^{r}=\sum_{l=0}^{n}{n+r \brace l+r}_{r,\lambda}(a^{+})^{l+r}a^{l},
\end{displaymath}
where ${n+r \brace l+r}_{r,\lambda}$ are the degenerate $r$-Stirling numbers of the second kind (see \eqref{13}).
In Theorem 3, we find that the differential operator $x^{k}D^{k}$ applied to the formal power series $f(x)$ can be represented as a linear combination of $(xD+r)_{m,\lambda}f(x)$ with the coefficients given by the unsigned degenerate $r$-Stirling numbers of the first kind. Hence we get the inverse relation to the normal ordering of $(a^{+}a+r)_{n,\lambda}$ (see Theorem 4). The inverse relation to Theorem 3 is derived in Theorem 5 in connection with the degenerate $r$-Stirling numbers of the second kind and hence we have the normal ordering of $(a^{+}a+r)_{n,\lambda}$  (see Theorem 6). In Theorems 7 and 9, two recurrence relations are obtained for the degenerate $r$-Stirling numbers of the second kind. Finally, we get a recurrence relation for the degenerate $r$-Stirling numbers of the first kind. For the rest of this section, we recall the facts that are needed throughout this paper.\par
It is well known that the Stirling numbers of the first kind are defined by 
\begin{equation} 
(x)_{n}=\sum_{k=0}^{n}S_{1}(n,k)x^{k},\quad (n\ge 0),\quad (\mathrm{see}\ [4-17]),	\label{1}
\end{equation}
where $(x)_{0}=1,\ (x)_{n}=x(x-1)\cdots(x-n+1),\ (n\ge 1)$. \par 
As the inversion formula of \eqref{1}, the Stirling numbers of the second kind are given by 
\begin{equation}
x^{n}=\sum_{k=0}^{n}S_{2}(n,k)(x)_{k},\quad (n\ge 0),\quad (\mathrm{see}\ [4,5,17]). \label{2}
\end{equation}
The unsigned Stirling number of the first kind is defined by 
\begin{displaymath}
	{n \brack k}=(-1)^{n-k}S_{1}(n,k),\quad (n,k\ge 0).
\end{displaymath}
For any $\lambda\in\mathbb{R}$, the degenerate factorial sequence is defined by 
\begin{equation}
(x)_{0,\lambda}=1,\quad (x)_{n,\lambda}=x(x-\lambda)	(x-2\lambda)\cdots(x-(n-1)\lambda),\quad (n\ge 1),\quad (\mathrm{see}\ [8,10]).\label{3}
\end{equation}
In [8], the degenerate Stirling numbers of the first kind are defined by 
\begin{equation}
(x)_{n}=\sum_{k=0}^{n}S_{1,\lambda}(n,k)(x)_{k,\lambda},\quad (n\ge 0). \label{4}	
\end{equation}
In view of \eqref{2}, the degenerate Stirling numbers of the second kind are given by 
\begin{equation}
(x)_{n,\lambda}=\sum_{k=0}^{n}S_{2,\lambda}(n,k)(x)_{k},\quad (n\ge 0),\quad (\mathrm{see}\ [8]).\label{5}
\end{equation}
The unsigned degenerate Stirling numbers of the first kind are given by 
\begin{equation}
\langle x\rangle_{n}=\sum_{k=0}^{n}{n \brack k}_{\lambda}\langle x\rangle_{k,\lambda},\quad (n\ge 0),\quad (\mathrm{see}\ [13-15]),\label{6}
\end{equation}
where 
\begin{displaymath}
	\langle x\rangle_{n}=(-1)^{n}(-x)_{n}\quad\mathrm{and}\quad \langle x\rangle_{n,\lambda}=(-1)^{n}(-x)_{n,\lambda},\quad (n\ge 0). 
\end{displaymath}
Note that 
\begin{displaymath}
	{n \brack k}_{\lambda}=(-1)^{n-k}S_{1,\lambda}(n,k).
\end{displaymath}
We note that the boson annihilation and creation operators $a$ and $a^{+}$ satisfy the commutation relation given by
\begin{equation}
[a,a^{+}]=aa^{+}-a^{+}a=1,\quad (\mathrm{see}\ [1,6,9,11]).\label{7}	
\end{equation}
In the corresponding Fock space, we have 
\begin{equation}
a^{+}|k\rangle=\sqrt{k+1}|k+1\rangle,\quad a|k\rangle=\sqrt{k}|k-1\rangle,\quad (\mathrm{see}\ [1,12]).\label{8}	
\end{equation}
Now, the number operator $\hat{n}$ is defined by 
\begin{equation}
\hat{n}|k\rangle=k|k\rangle.\label{9}
\end{equation}
By \eqref{8}, we get $a^{+}a|k\rangle=k|k\rangle$, that is $\hat{n}=a^{+}a$. Thus, we note that 
\begin{displaymath}
	[a,\hat{n}]=a\hat{n}-\hat{n}a=a,\quad [\hat{n},a^{+}]=\hat{n}a^{+}-a^{+}\hat{n}=a^{+}. 
\end{displaymath}
The normal ordering of an integral power of the number operator $\hat{n}=a^{+}a$ in terms of boson operators $a$ and $a^{+}$ can be written in the form 
\begin{equation}
(a^{+}a)^{k}=\sum_{l=0}^{k}S_{2}(k,l)(a^{+})^{l}a^{l},\quad (k\ge 0),\quad (\mathrm{see}\ [1,6,16]).\label{9-1}
\end{equation}
Recently, the normal ordering of the degenerate $k$-th power of the number operator $a^{+}a$, namely $(a^{+}a)_{k,\lambda}$, in terms of boson operators $a$ and $a^{+}$ can be written in the form 
\begin{equation}
(a^{+}a)_{k,\lambda}=\sum_{l=0}^{k}S_{2,\lambda}(k,l)(a^{+})^{l}a^{l},\quad (k\ge 0),\quad (\mathrm{see}\ [12]).\label{9-2}	
\end{equation}
For $r\ge 0$, the unsigned $r$-Stirling numbers are defined by 
\begin{equation}
\langle x+r\rangle_{n}=\sum_{k=0}^{n}{n+r \brack k+r}_{r}x^{k},\quad (n\ge 0),\quad (\mathrm{see}\ [13-15]),\label{10}	
\end{equation}
where 
\begin{displaymath}
	\langle x\rangle_{n}=x(x+1)\cdots(x+n-1),\quad (n\ge 1),\quad \langle x\rangle_{0}=1.
\end{displaymath}
In view of \eqref{2}, the $r$-Stirling numbers of the second kind are given by 
\begin{equation}
(x+r)^{n}=\sum_{k=0}^{n}{n+r \brace k+r}_{r}(x)_{k},\quad (n\ge 0),\quad (\mathrm{see}\ [14,15]).\label{11}
\end{equation}
Recently, the authors introduced the unsigned degenerate $r$-Stirling numbers of the first kind, which are given by 
\begin{equation}
\langle x+r\rangle_{n}=\sum_{k=0}^{n}{n+r \brack k+r}_{r,\lambda}\langle x\rangle_{k,\lambda},\quad (n\ge 0),\quad (\mathrm{see}\ [13-15]),\label{12}
\end{equation}
In view of \eqref{5}, the degenerate $r$-Stirling numbers of the second kind are defined by the authors, which are given by 
\begin{equation}
(x+r)_{n,\lambda}=\sum_{k=0}^{n}{n+r \brace k+r}_{r,\lambda}(x)_{k},\quad (n\ge 0),\quad (\mathrm{see}\ [13-15]). \label{13}	
\end{equation}
In addition, the degenerate $r$-Stirling numbers of the first kind are given by 
\begin{equation}
(x-r)_{n}=\sum_{k=0}^{n}S_{1,\lambda}^{(r)}(n+r,k+r)(x)_{k,\lambda},\quad (n\ge 0),\quad (\mathrm{see}\ [13,15]).\label{14}	
\end{equation}
Note that 
\begin{displaymath}
	{n+r \brack k+r}_{r,\lambda}=(-1)^{n-k}S_{1,\lambda}^{(r)}(n+r,k+r),\quad (n\ge 0). 
\end{displaymath}
In this paper, we study the degenerate $r$-Stirling numbers associated with the number operator $\hat{n}=a^{+}a$ in terms of boson operators $a$ and $a^{+}$ that satisfy $[a,a^{+}]=aa^{+}-a^{+}a=1$. \par 
That is, we derive some identities involving the degenerate $r$-Stirling numbers arising from the normal ordering of degenerate integral powers of number operator. $\hat{n}=a^{+}a$ in terms of boson operators $a$ and $a^{+}$, which are given by 
\begin{displaymath}
	(a^{+}a+r)_{n,\lambda}=\sum_{k=0}^{n}{n+r \brace k+r}_{r,\lambda}(a^{+})^{k}a^{k}. 
\end{displaymath}

\section{Some identities on degenerate $r$-Stirling numbers via boson operators}
Let $f(x)$ be a formal power series in $x$ with coefficients in $\mathbb{C}$, given by $f(x)=\sum_{n=0}^{\infty}a_{n}x^{n}$. For $r\ge 0$, we observe that 
\begin{align}
(xD)_{m,\lambda}x^{r}f(x)&=(xD)(xD-\lambda)\cdots (xD-(m-1)\lambda)x^{r}\sum_{k=0}^{\infty}a_{k}x^{k}\label{15}\\
&=\sum_{k=0}^{\infty}a_{k}(k+r)_{m,\lambda}x^{k+r}=\sum_{k=0}^{\infty}a_{k}\sum_{l=0}^{m}{m+r \brace l+r}_{r,\lambda}(k)_{l}x^{k+r} \nonumber \\
&=\sum_{l=0}^{m}{m+r \brace l+r}_{r,\lambda}\sum_{k=0}^{\infty}a_{k}(k)_{l}x^{k}x^{r}=\sum_{l=0}^{m}{m+r \brace l+r}_{r,\lambda}x^{l+r}f^{(l)}(x), \nonumber
\end{align}
where 
\begin{displaymath}
f^{(l)}(x)=\bigg(\frac{d}{dx}\bigg)^{l}f(x),\quad D=\frac{d}{dx}.
\end{displaymath}
{\it{From now on throughout this paper, $D$ will denote the operator $D=\frac{d}{dx}$.}} \par
Now, by \eqref{15}, we get 
\begin{equation}
(xD)_{m,\lambda}x^{r}f(x)=\sum_{l=0}^{m}{m+r \brace l+r}_{r,\lambda}\Big(D^{l}f(x)\Big)x^{l+r}.\label{16}
\end{equation}
In particular, for $m,r\ge 0$ with $m\ge r$, we have 
\begin{equation}
(xD)_{m-r,\lambda}x^{r}f(x)=\sum_{l=r}^{m}{m \brace l}_{r,\lambda}\Big( D^{l-r}f(x)\Big)x^{l}.\label{17}
\end{equation}
Replacing $f(x)$ by $f^{(r)}(x)$, we get 
\begin{equation}
(xD)_{m-r,\lambda}x^{r}D^{r}f(x)=\sum_{l=r}^{m}{m \brace l}_{r,\lambda}x^{l}D^{l}f(x). \label{18}
\end{equation}
Therefore, by \eqref{16} and \eqref{18}, we obtain the following theorem. 
\begin{theorem}
	For $m,r\in\mathbb{Z}$ with $m, r \ge 0$, we have 
	\begin{displaymath}
		(xD)_{m,\lambda}x^{r}f(x)=\sum_{l=0}^{m}{m+r \brace l+r}_{r,\lambda}x^{l+r}D^{l}f(x). 
	\end{displaymath}
	In addition, for $m,r\in\mathbb{Z}$ with $m\ge r\ge 0$, we have 
	\begin{displaymath}
		(xD)_{m-r,\lambda}x^{r}D^{r}f(x)=\sum_{l=r}^{m}{m \brace l}_{r,\lambda}x^{l}D^{l}f(x).
	\end{displaymath}
\end{theorem}
Note that $Dx-xD=1$. We recall that the bosonic commutation relation $[a,a^{+}]=aa^{+}-a^{+}a=1$ can be realized formally in a suitable space of functions by letting $a=\frac{d}{dx}$ and $a^{+}=x$.  
Therefore we obtain the following corollary. 
\begin{corollary}
For $m,r\in\mathbb{Z}$ with $m, r \ge 0$, we have 
\begin{displaymath}
(a^{+}a)_{m,\lambda}(a^{+})^{r}=\sum_{l=0}^{m}{m+r \brace l+r}_{r,\lambda}(a^{+})^{l+r}a^{l}. 
\end{displaymath}
In addition, for $m,r\in\mathbb{Z}$ with $m\ge r\ge 0$, we have 
\begin{displaymath}
(a^{+}a)_{m-r,\lambda}(a^{+})^{r}a^{r}=\sum_{l=r}^{m}{m\brace l}_{r,\lambda}(a^{+})^{l}a^{l}. 
\end{displaymath}
\end{corollary}
By \eqref{13}, we get 
\begin{align}
\sum_{k=0}^{n+1}{n+r+1 \brace k+r}_{r,\lambda}(x)_{k}&=(x+r)_{n+1,\lambda}=(x+r)_{n,\lambda}(x+r-n\lambda) \label{20} \\
&=\sum_{k=0}^{n}{n+r \brace k+r}_{r,\lambda}(x)_{k}(x-k+k+r-n\lambda) \nonumber \\
&=\sum_{k=0}^{n}{n+r \brace k+r}_{r,\lambda}(x)_{k+1}+\sum_{k=0}^{n}{n+r \brace k+r}_{r,\lambda}(k+r-n\lambda)(x)_{k} \nonumber \\
&=\sum_{k=0}^{n+1}\bigg({n+r \brace k+r-1}_{r,\lambda}+(k+r-n\lambda){n+r \brace k+r}_{r,\lambda}\bigg)(x)_{k}.\nonumber	
\end{align}
Comparing the coefficients on both sides of \eqref{20}, we have 
\begin{displaymath}
	{n+r+1 \brace k+r}_{r,\lambda}={n+r \brace k+r-1}_{r,\lambda}+(k+r-n\lambda){n+r \brace k+r}_{r,\lambda},
\end{displaymath}
where $n,k\in\mathbb{Z}$ with $n+1\ge k\ge 0$. \par 
Now, we observe that 
\begin{align}
x^{k}D^{k}f(x)&=\sum_{n=0}^{\infty}a_{n}(n)_{k}x^{n}\label{21} \\
&=\sum_{n=0}^{\infty}a_{n}\bigg(\sum_{m=0}^{k}S_{1,\lambda}^{(r)}(k+r,m+r)(n+r)_{m,\lambda}\bigg)x^{n} \nonumber \\
&=\sum_{m=0}^{k}S_{1,\lambda}^{(r)}(k+r,m+r)\sum_{n=0}^{\infty}a_{n}(n+r)_{m,\lambda}x^{n} \nonumber\\
&=\sum_{m=0}^{k}S_{1,\lambda}^{(r)}(k+r,m+r)\bigg(x\frac{d}{dx}+r\bigg)_{m,\lambda}f(x) \nonumber \\
&=\sum_{m=0}^{k}(-1)^{k-m}{k+r \brack m+r}_{r,\lambda}\bigg(x\frac{d}{dx}+r\bigg)_{m,\lambda}f(x). \nonumber
\end{align}
Therefore, by \eqref{21}, we obtain the following theorem. 
\begin{theorem}
	For $k\ge 0$, we have 
	\begin{equation}
	x^{k}D^{k}f(x)=\sum_{m=0}^{k}(-1)^{k-m}{k+r \brack m+r}_{r,\lambda}(xD+r)_{m,\lambda}f(x).\label{22}	\end{equation}
\end{theorem}
In \eqref{22}, the normal ordering of degenerate integral power of the number operator $\hat{n}a^{+}a$ in terms of boson operators $a^{+}$ and $a$ can be presented in the form 
\begin{align}
(a^{+})^{k}a^{k}&=\sum_{m=0}^{k}(-1)^{k-m}{k+r \brack m+r}_{r,\lambda}(a^{+}a+r)_{m,\lambda}\label{23} \\
&=\sum_{m=0}^{k}(-1)^{k-m}{k+r \brack m+r}_{r,\lambda}(\hat{n}+r)_{m,\lambda}\nonumber.
\end{align}
For $k,r\in\mathbb{Z}$ with $k\ge r\ge 0$, by \eqref{23}, we get 
\begin{align}
(a^{+})^{k-r}a^{k-r}&=\sum_{m=r}^{k}(-1)^{k-m}{k \brack m}_{r,\lambda}(a^{+}a+r)_{m-r,\lambda} \label{24} \\
&=\sum_{m=r}^{k}(-1)^{k-m}{k \brack m}_{r,\lambda}(\hat{n}+r)_{m-r,\lambda}. \nonumber	
\end{align}
Therefore, by \eqref{23} and \eqref{24}, we obtain the following theorem. 
\begin{theorem}
For $k,r\in\mathbb{Z}$ with $k\ge r\ge 0$, we have 
\begin{displaymath}
	(a^{+})^{k}a^{k}=\sum_{m=0}^{k}(-1)^{k-m}{k+r \brack m+r}_{r,\lambda}(\hat{n}+r)_{m,\lambda},
\end{displaymath}
and
\begin{displaymath}
	(a^{+})^{k-r}a^{k-r}=\sum_{m=r}^{k}(-1)^{k-m}{k \brack m}_{r,\lambda}(\hat{n}+r)_{m-r,\lambda}.
\end{displaymath}
\end{theorem}
As the inversion formula of \eqref{22}, we consider the following differential equation:
\begin{align}
\bigg(x\frac{d}{dx}+r\bigg)_{n,\lambda}f(x)&=\sum_{m=0}^{\infty}a_{m}\bigg(x\frac{d}{dx}+r\bigg)_{n,\lambda}x^{m}=\sum_{m=0}^{\infty}a_{m}(m+r)_{n,\lambda}x^{m}\label{25} \\
&=\sum_{m=0}^{\infty}a_{m}\sum_{k=0}^{n}{n+r \brace k+r}_{r,\lambda}(m)_{k}x^{m} \nonumber \\
&=\sum_{k=0}^{n}{n+r \brace k+r}_{r,\lambda}\sum_{m=0}^{\infty}a_{m}(m)_{k}x^{m} \nonumber \\
&=\sum_{k=0}^{n}{n+r \brace k+r}_{r,\lambda}x^{k}\bigg(\frac{d}{dx}\bigg)^{k}f(x).\nonumber
\end{align}
Therefore, by \eqref{25}, we obtain the following theorem. 
\begin{theorem}
	For $n\ge 0$, we have 
	\begin{equation}
	(xD+r)_{n,\lambda}f(x)=\sum_{k=0}^{n}{n+r \brace k+r}_{r,\lambda}x^{k}D^{k}f(x). \label{26}
	\end{equation}
\end{theorem}
From Theorems 1 and 5, we note that 
\begin{equation}
x^{r}(xD+r)_{n,\lambda}=(xD)_{n,\lambda}x^{r},\quad (n\ge 0). \label{27}	
\end{equation}
Thus, by \eqref{27}, we get 
\begin{equation}
(a^{+})^{r}(a^{+}a+r)_{n,\lambda}=(a^{+}a)_{n,\lambda}(a^{+})^{r},\quad (n\ge 0). \label{28}
\end{equation}
From \eqref{26}, we have 
\begin{equation}
(a^{+}a+r)_{n,\lambda}=\sum_{k=0}^{n}{n+r \brace k+r}_{r,\lambda}(a^{+})^{k}a^{k}. \label{29}
\end{equation}
Therefore, by \eqref{28} and \eqref{29}, we obtain the following theorem. 
\begin{theorem}
	For $n\ge 0$, we have 
	\begin{displaymath}
		(a^{+})^{r}(a^{+}a+r)_{n,\lambda}=(a^{+}a)_{n,\lambda}(a^{+})^{r}, 
	\end{displaymath}
	and 
	\begin{equation}
		(a^{+}a+r)_{n,\lambda}=\sum_{k=0}^{n}{n+r \brace k+r}_{r,\lambda}(a^{+})^{k}a^{k}. \label{30}
	\end{equation}
\end{theorem}
By \eqref{30}, we get 
\begin{align}
a^{+}(\hat{n}+r+1-\lambda)_{k,\lambda}a&=a^{+}\sum_{l=0}^{k}\binom{k}{l}(\hat{n}+r)_{l,\lambda}(1-\lambda)_{k-l,\lambda}a \label{31} \\
&=\sum_{l=0}^{k}\binom{k}{l}(1-\lambda)_{k-l,\lambda}\sum_{m=0}^{l}{l+r \brace m+r}_{r,\lambda}(a^{+})^{m+1}a^{m+1} \nonumber \\
&=\sum_{m=0}^{k}\sum_{l=m}^{k}\binom{k}{l}(1-\lambda)_{k-l,\lambda}{l+r \brace m+r}_{r,\lambda}(a^{+})^{m+1}a^{m+1}. \nonumber
\end{align}
On the other hand, by \eqref{29} and \eqref{31}, we get 
\begin{align}
&\sum_{m=0}^{k+1}{k+1+r \brace m+r}_{r,\lambda}(a^{+})^{m}a^{m}=(\hat{n}+r)_{k+1,\lambda}=(\hat{n}+r)(\hat{n}+r-\lambda)_{k,\lambda} \label{32} \\
&=\hat{n}(\hat{n}+r-\lambda)_{k,\lambda}+r(\hat{n}+r-\lambda)_{k,\lambda}=a^{+}(\hat{n}+r+1-\lambda)_{k,\lambda}a+r(\hat{n}+r-\lambda)_{k,\lambda}\nonumber\\
&=\sum_{m=0}^{k}\sum_{l=m}^{k}\binom{k}{l}(1-\lambda)_{k-l,\lambda}{l+r \brace m+r}_{r,\lambda}(a^{+})^{m+1}a^{m+1}+r\sum_{l=0}^{k}\binom{k}{l}(-\lambda)_{k-l,\lambda}(\hat{n}+r)_{l,\lambda}\nonumber \\
&=\sum_{m=1}^{k+1}\sum_{l=m-1}^{k}\binom{k}{l}(1-\lambda)_{k-l,\lambda}{l+r \brace m-1+r}_{r,\lambda}(a^{+})^{m}a^{m}+r\sum_{l=0}^{k}\binom{k}{l}(-\lambda)_{k-l,\lambda}\sum_{m=0}^{l}{l+r \brace m+r}_{r,\lambda}(a^{+})^{m}a^{m}\nonumber \\
&=\sum_{m=1}^{k+1}\sum_{l=m-1}^{k}\binom{k}{l}(1-\lambda)_{k-l,\lambda}{l+r \brace m+r-1}_{r,\lambda}(a^{+})^{m}a^{m}+\sum_{m=0}^{k}\bigg(r\sum_{l=m}^{k}\binom{k}{l}(-\lambda)_{k-l,\lambda}{l+r \brace m+r}_{r,\lambda}\bigg)(a^{+})^{m}a^{m} \nonumber \\
&=\sum_{m=0}^{k+1}\bigg(\sum_{l=m-1}^{k}\binom{k}{l}(1-\lambda)_{k-l,\lambda}{l+r \brace m+r-1}_{r,\lambda}+r\sum_{l=m}^{k}\binom{k}{l}(-\lambda)_{k-l,\lambda}{l+r \brace m+r}_{r,\lambda}\bigg)(a^{+})^{m}a^{m}. \nonumber
\end{align}
Therefore, by comparing the coefficients on both sides of \eqref{32}, we obtain the following theorem. 
\begin{theorem}
	For $m,k\in\mathbb{Z}$ with $k\ge m\ge 0$, we have 
	\begin{displaymath}
		{k+1+r \brace m+r}_{r,\lambda}=\sum_{l=m-1}^{k}\binom{k}{l}(1-\lambda)_{k-l,\lambda}{l+r \brace m+r-1}_{r,\lambda}+r\sum_{l=m}^{k}\binom{k}{l}(-\lambda)_{k-l,\lambda}{l+r \brace m+r}_{r,\lambda}.
	\end{displaymath}
\end{theorem}
We observe that 
\begin{align}
\hat{n}(\hat{n}+r-1)_{k,\lambda}&=a^{+}a(a^{+}a+r-1)_{k,\lambda}\label{33}\\
&=a^{+}a(a^{+}a+r-1)(a^{+}a+r-1-\lambda)\cdots(a^{+}a+r-1-(k-1)\lambda)\nonumber \\
&=a^{+}(aa^{+}+r-1)a(a^{+}a+r-1-\lambda)(a^{+}a+r-1-2\lambda)\cdots(a^{+}a+r-1-(k-1)\lambda)\nonumber \\
&=a^{+}(aa^{+}+r-1)(aa^{+}+r-1-\lambda)a(a^{+}a+r-1-2\lambda)\cdots(a^{+}a+r-1-(k-1)\lambda)\nonumber \\
&=\cdots\nonumber \\
&=a^{+}(aa^{+}+r-1)(aa^{+}+r-1-\lambda)(aa^{+}+r-1-2\lambda)\cdots(aa^{+}+r-1-(k-1)\lambda)a\nonumber \\
&=a^{+}(a^{+}a+r)(a^{+}a+r-\lambda)(a^{+}a+r-2\lambda)\cdots(a^{+}a+r-(k-1)\lambda)a\nonumber \\
&=a^{+}(a^{+}a+r)_{k,\lambda}a=a^{+}(\hat{n}+r)_{k,\lambda}a. \nonumber
\end{align}
Therefore, by \eqref{33}, we obtain the following lemma. 
\begin{lemma}
	For $k\ge 0$, we have 
	\begin{equation}
	\hat{n}(\hat{n}+r-1)_{k,\lambda}=a^{+}(\hat{n}+r)_{k,\lambda}a.\label{34}	
	\end{equation}
\end{lemma}
From \eqref{34} and Theorem 6, we note that 
\begin{align}
&a^{+}(\hat{n}+r)_{k,\lambda}a=\hat{n}(\hat{n}+r-1)_{k,\lambda}=\sum_{m=0}^{k}\binom{k}{m}(-1)_{k-m,\lambda}\hat{n}(\hat{n}+r)_{m,\lambda}	\label{35} \\
&=\sum_{m=0}^{k}\binom{k}{m}(-1)^{k-m}\langle 1\rangle_{k-m,\lambda}(\hat{n}+r-m\lambda+m\lambda-r)(\hat{n}+r)_{m,\lambda}\nonumber \\
&=\sum_{m=0}^{k}\binom{k}{m}(-1)^{k-m}\langle 1\rangle_{k-m,\lambda}\Big((\hat{n}+r)_{m+1,\lambda}+(m\lambda-r)(\hat{n}+r)_{m,\lambda}\Big) \nonumber \\
&=\sum_{m=0}^{k}\binom{k}{m}(-1)^{k-m}\langle 1\rangle_{k-m,\lambda}(\hat{n}+r)_{m+1,\lambda}\nonumber \\
&\qquad +\lambda\sum_{m=0}^{k}\binom{k}{m+1}(-1)^{k-m-1}\langle 1\rangle_{k-m-1,\lambda}(m+1)(\hat{n}+r)_{m+1,\lambda}\nonumber \\
&\qquad -r\sum_{m=0}^{k}\binom{k}{m}(-1)^{k-m}\langle 1\rangle_{k-m,\lambda}(\hat{n}+r)_{m,\lambda}\nonumber \\
&=\sum_{m=0}^{k}\binom{k}{m}(-1)^{k-m}\langle 1\rangle_{k-m,\lambda}\sum_{p=0}^{m+1}{m+1+r \brace p+r}_{r,\lambda}(a^{+})^{p}a^{p} \nonumber \\
&\qquad +\lambda\sum_{m=0}^{k}\binom{k}{m+1}(-1)^{k-m-1}\langle 1\rangle_{k-m-1,\lambda}(m+1)\sum_{p=0}^{m+1}{m+1+r \brace p+r}_{r,\lambda}(a^{+})^{p}a^{p}\nonumber \\
&\qquad -r\sum_{m=0}^{k}\binom{k}{m}(-1)^{k-m}\langle 1\rangle_{k-m,\lambda}\sum_{p=0}^{m+1}{m+r \brace p+r}_{r,\lambda}(a^{+})^{p}a^{p} \nonumber \\
&=\sum_{p=1}^{k+1}\sum_{m=p-1}^{k}\binom{k}{m}(-1)^{k-m}\langle 1\rangle_{k-m,\lambda}{m+1+r \brace p+r}_{r,\lambda}(a^{+})^{p}a^{p} \nonumber \\
&\qquad +\lambda\sum_{p=1}^{k+1}\sum_{m=p-1}^{k}\binom{k}{m+1}(-1)^{k-m-1}\langle 1\rangle_{k-m-1,\lambda}(m+1){m+1+r \brace p+r}_{r,\lambda}(a^{+})^{p}a^{p}\nonumber 
\end{align}
\begin{align*}
&\qquad -r \sum_{p=1}^{r+1}\sum_{m=p-1}^{k}\binom{k}{m}(-1)^{k-m} \langle 1 \rangle_{k-m,\lambda}{ m+r \brace p+r}_{r,\lambda}(a^{+})^{p}a^{p} \nonumber \\
&=\sum_{p=0}^{k}\bigg\{\sum_{m=p}^{k}\binom{k}{m}(-1)^{k-m}\langle 1\rangle_{k-m,\lambda}{m+1+r \brace p+1+r}_{r,\lambda}\nonumber\\
&\qquad +\lambda\sum_{m=p}^{k}\binom{k}{m+1}(-1)^{k-m-1}\langle 1\rangle_{k-m-1}(m+1){m+1+r \brace p+1+r}_{r,\lambda}\nonumber \\
&\qquad -r\sum_{m=p}^{k}\binom{k}{m}(-1)^{k-m}\langle 1\rangle_{k-m,\lambda}{m+r \brace p+1+r}_{r,\lambda}\bigg\}(a^{+})^{p+1}a^{p+1}. \nonumber
\end{align*}
On the other hand, by \eqref{29}, we get 
\begin{equation}
a^{+}(\hat{n}+r)_{k,\lambda}a=\sum_{p=0}^{k}{k+r \brace p+r}_{r,\lambda}(a^{+})^{p+1}a^{p+1}. \label{36}	
\end{equation}
Therefore, by \eqref{35} and \eqref{36}, we obtain the following theorem. 
\begin{theorem}
For $p,k\in\mathbb{Z}$ with $k\ge p\ge 0$, we have 
\begin{align*}
&{k+r \brace p+r}_{r,\lambda}=\sum_{m=p}^{k}\bigg[\binom{k}{m}(-1)^{k-m}\langle 1\rangle_{k-m,\lambda}\bigg({m+1+r\brace p+1+r}_{r,\lambda}-r{m+r \brace p+1+r}_{r,\lambda}\bigg)\\
&\qquad \qquad \qquad\qquad +\lambda k\binom{k-1}{m}\langle 1\rangle_{k-m-1}(-1)^{k-m-1}{m+1+1 \brace p+1+r}_{r,\lambda}\bigg].
\end{align*}
\end{theorem}
From Theorem 4, we note that 
\begin{align}
(a^{+})^{k+1}a^{k+1}&=\sum_{l=0}^{k}(-1)^{k-l}{k+r \brack l+r}_{r,\lambda}a^{+}(\hat{n}+r)_{l,\lambda}a\label{37}\\
&=\sum_{l=0}^{k}(-1)^{k-l}{k+r \brack l+r}_{r,\lambda}\hat{n}(\hat{n}+r-1)_{l,\lambda}\nonumber \\
&=\sum_{l=0}^{k}{k+r \brack l+r}_{r,\lambda}(-1)^{k-l}\sum_{m=0}^{l}\binom{l}{m}(-1)^{l-m}\langle 1\rangle_{l-m,\lambda}(\hat{n}+r)_{m,\lambda}\hat{n}\nonumber\\
&=\sum_{l=0}^{k}{k+r \brack l+r}_{r,\lambda}(-1)^{k-l}\sum_{m=0}^{l}\binom{l}{m}(-1)^{l-m}\langle 1\rangle_{l-m,\lambda}(\hat{n}+r)_{m,\lambda}(\hat{n}+r-m\lambda+m\lambda-r)\nonumber \\
&=\sum_{l=0}^{k}{k+r \brack l+r}_{r,\lambda}(-1)^{k-l}\sum_{m=0}^{l}\binom{l}{m}(-1)^{l-m}\langle 1\rangle_{l-m,\lambda}(\hat{n}+r)_{m+1,\lambda}\nonumber \\
&\quad +\sum_{l=0}^{k}{k+r \brack l+r}_{r,\lambda}(-1)^{k-l}\sum_{m=0}^{l}\binom{l}{m}(-1)^{l-m}\langle 1\rangle_{l-m,\lambda}(\hat{n}+r)_{m,\lambda}(m\lambda-r)\nonumber \\
&=\sum_{m=0}^{k}\sum_{l=m}^{k}{k+r \brack l+r}_{r,\lambda}(-1)^{k-m}\binom{l}{m}\langle 1\rangle_{l-m,\lambda}(\hat{n}+r)_{m+1,\lambda}\nonumber
\end{align}
\begin{align*}
&\quad +\sum_{m=0}^{k}\sum_{l=m}^{k}{k+r \brack l+r}_{r,\lambda}(-1)^{k-m}\binom{l}{m}\langle 1\rangle_{l-m,\lambda}(m\lambda-r)(\hat{n}+r)_{m,\lambda}\nonumber \\
&=\sum_{m=1}^{k+1}\bigg(\sum_{l=m-1}^{k}{k+r \brack l+r}_{r,\lambda}(-1)^{k-m+1}\binom{l}{m-1}\langle 1\rangle_{l-m+1,\lambda}\bigg)(\hat{n}+r)_{m,\lambda}\nonumber \\
&\quad +\sum_{m=0}^{k}(m\lambda-r)\sum_{l=m}^{k}{k+r \brack l+r}_{r,\lambda}(-1)^{k-m}\binom{l}{m}\langle 1\rangle_{l-m,\lambda}(\hat{n}+r)_{m,\lambda}\nonumber \\
&=\sum_{m=0}^{k+1}\bigg\{\sum_{l=m-1}^{k}{k+r \brack l+r}_{r,\lambda}(-1)^{k-m+1}\binom{l}{m-1}\langle 1\rangle_{l-m+1,\lambda}\nonumber \\
&\qquad +(m\lambda-r)\sum_{l=m}^{k}{k+r \brack l+r}_{r,\lambda}(-1)^{k-m}\binom{l}{m}\langle 1\rangle_{l-m,\lambda}\bigg\}(\hat{n}+r)_{m,\lambda}.\nonumber
\end{align*}
On the other hand, by \eqref{23}, we get 
\begin{align}
(a^{+})^{k+1}a^{k+1}&=\sum_{m=0}^{k+1}(-1)^{k-m+1}{k+r+1 \brack m+r}_{r,\lambda}	(a^{+}a+r)_{m,\lambda}\label{38}\\
&=\sum_{m=0}^{k+1}(-1)^{k-m+1}{k+r+1 \brack m+r}_{r,\lambda}(\hat{n}+r)_{m,\lambda}.\nonumber
\end{align}
Therefore, by \eqref{37} and \eqref{38}, we obtain the following theorem. 
\begin{theorem}
	For $m,k\in\mathbb{Z}$ with $k\ge m\ge 0$, we have 
	\begin{align*}
		(-1)^{k-m+1}{k+r+1 \brack m+r}_{r,\lambda}&=\sum_{l=m-1}^{k}{k+r \brack l+r}_{r,\lambda}(-1)^{k-m+1}\binom{l}{m-1}\langle 1\rangle_{l-m+1,\lambda}\nonumber \\
		&\quad +(m\lambda-r)\sum_{l=m}^{k}{k+r \brack l+r}_{r,\lambda}(-1)^{k-m}\binom{l}{m}\langle 1\rangle_{l-m,\lambda}.
	\end{align*}
\end{theorem}

\section{conclusion}
We have witnessed that various degenerate versions of some special numbers and polynomials were studied by means of many different methods, which include combinatorial methods, generating functions, $p$-adic analysis, umbral calculus, probability theory, differential equations, operator theory and analytic number theory.\par
The Stirling numbers of the second kind appear as the coefficients in the normal ordering of an integral power of the number operator in terms of the boson operators. Recently, it was shown that the degenrate Stirling numbers of the second kind appear as the coefficients in the normal ordering of a degenerate integral power of the number operator in terms of boson operators. In this paper, we obtained the normal ordering of a degenerate integral power of the number operator multiplied by an integral power of the creation boson operator in trems of boson operators where the degenerate $r$-Stirling numbers of the second kind appear as the coefficients. In addition, we derived recurrence relations for the degenerate $r$-Stirling numbers of both kinds from certain normal orderings for boson operators.\par
It is one of our future research projects to continue to explore various degenerate versions of some special numbers and polynomials and to find their applications in physics, science and engineering as well as in mathematics.

\end{document}